\documentclass[12pt,twoside]{article}

\usepackage{amssymb, amsthm, amsmath}
\newtheorem{theorem}{Theorem}
\newtheorem{definition}[theorem]{Definition}

\newtheorem{corollary}[theorem]{Corollary}

\begin{document}
\centerline {\Large{\bf The generalized r-Whitney numbers }}
\centerline{}
\centerline{\bf {B. S. El-Desouky}}

\centerline{}

\centerline{Department of Mathematics, Faculty of Science}

\centerline{Mansoura University, Mansoura 35516, Egypt}

\centerline{b\_desouky@yahoo.com}
\centerline{}

\centerline{\bf {F. A. Shiha}}

\centerline{}

\centerline{Department of Mathematics, Faculty of Science}

\centerline{Mansoura University, Mansoura 35516, Egypt}

\centerline{fshiha@yahoo.com}
\centerline{}

\centerline{\bf {Ethar M. Shokr}}

\centerline{}

\centerline{Department of Mathematics, Faculty of Science}

\centerline{Mansoura University, Mansoura 35516, Egypt}

\centerline{ethar.shokr@gmail.com}
\centerline{}
\begin{abstract} In this paper, we define the generalized r-Whitney numbers of the first and second kind. Moreover, we drive the generalized Whitney numbers of the first and second kind. The recurrence relations and the generating functions of these numbers are derived. The relations between these numbers and generalized Stirling numbers of the first and second kind are deduced. Furthermore, some special cases are given. Finally, matrix representation of The relations between Whitney and Stirling numbers are given.
\end{abstract}
{\bf Keywords:} Whitney numbers; r-Whitney numbers; p-Stirling numbers; generalized q-Stirling numbers; generalized Stirling numbers.\\
{\bf Mathematics Subject Classification:} 05A19, 11B73, 11B75, 11C20. \\
\section{Introduction}
\label{intro}
The Whitney numbers of the first and second kind were introduced, respectively, by Dowling \cite{Dow1} and Benoumhani \cite{Ben1} as
\[
m^n(x)_n=\sum_{k=0}^{n}w_{m,r}(n,k)(mx+1)^k,
\]
and
\[
(mx+1)^n=\sum_{k=0}^{n}W_{m,r}(n,k) m^k(x)_k.
\]
The r-Whitney numbers of the first and second kind were introduced, respectively, by Mez\H o \cite{Mez1} as
\[
m^n(x)_n=\sum_{k=0}^{n}w_{m,r}(n,k)(mx+r)^k,
\]
and
\[
(mx+r)^n=\sum_{k=0}^{n}W_{m,r}(n,k) m^k(x)_k.
\]
Let $\overline{\boldsymbol{\alpha}}=(\alpha_{0},\alpha_{1},\cdots,\alpha_{n-1})$ where $\alpha_i ,\:i=0,1,\cdots n-1$ are real numbers.\\
 In this paper we use the following notations
\begin{equation}
(x;\overline{\boldsymbol{\alpha}})_n=\prod_{i=0}^{n-1}(x-\alpha_i),\;(x;\overline{\boldsymbol{\alpha}})_0=1,
\end{equation}
\begin{eqnarray*}
[x;\overline{\boldsymbol{\alpha}}]_{\underline{n},q}&=&\prod_{i=0}^{n-1}[x-\alpha_i]_q\\
&=&q^{-\sum_{i=0}^{n-1}\alpha_i}([x]_q;\overline{[\boldsymbol{\alpha}]}_q)_n,
\end{eqnarray*}
where $([x]_q;\overline{[\boldsymbol{\alpha}]}_q)_n=([x]_q-[\alpha_0]_q)([x]_q-[\alpha_1]_q)\cdots([x]_q-[\alpha_{n-1}]_q)$, see \cite{Cha1}.\\
\\
and
\begin{equation}
(\alpha_i)_n=\prod_{j=0,j\neq{i}}^{n-1}(\alpha_i-\alpha_j),\;(\alpha_i)_0=1,
\end{equation}
see El-Desouky and Caki\'{c} \cite{Des1}, Comtet \cite{Com1}.\\
\\
This paper is organized as follows:\\
In Sections 2 and 3 we derive the generalized r-Whitney numbers of the first and second kind.The recurrence relations and the generating functions of these numbers are derived. Furthermore, some interesting special cases of these numbers are given.
In Section 4 we obtained the generalized Whitney numbers of the first and second kind by setting $r=1$.
We investigate some relations between the generalized r-Whitney numbers and Stirling numbers and generalized harmonic numbers in Section 5. Finally, we obtained a matrix representation for these relations.
\section{The generalized r-Whitney numbers of the first kind}
\label{sec:1}
\begin{definition} The generalized r-Whitney numbers of the first kind $w_{m,r;\overline{\boldsymbol{\alpha}}}(n,k)$ with parameter $\overline{\boldsymbol{\alpha}}=(\alpha_0, \alpha_1,\cdots, \alpha_{n-1})$ are defined by
\begin{equation}\label{1}
m^n(x;\overline{\boldsymbol{\alpha}})_n=\sum_{k=0}^{n}w_{m,r;\overline{\boldsymbol{\alpha}}}(n,k)(mx+r)^k,
\end{equation}
where $w_{m,r;\overline{\boldsymbol{\alpha}}}(0,0)=1$ and $w_{m,r;\overline{\boldsymbol{\alpha}}}(n,k)=0$ for $k>n$.
\end{definition}
\begin{theorem}
The generalized r-Whitney numbers of the first kind $w_{m,r;\overline{\boldsymbol{\alpha}}}(n,k)$ satisfy the recurrence relation
\begin{equation}\label{2}
w_{m,r;\overline{\boldsymbol{\alpha}}}(n,k)=w_{m,r;\overline{\boldsymbol{\alpha}}}(n-1,k-1)-(r+m\alpha_{n-1})w_{m,r;\overline{\boldsymbol{\alpha}}}(n-1,k),
\end{equation}
for $k\geq1$ and $w_{m,r;\overline{\boldsymbol{\alpha}}}(n,0)=(-1)^n\prod_{i=1}^{n}(r+m\alpha_{i-1})$.
\end{theorem}
{ \it Proof.}
Since $m^n(x;\overline{\boldsymbol{\alpha}})_n=(mx+r-r-m\alpha_{n-1})m^{n-1}(x;\overline{\boldsymbol{\alpha}})_{n-1}$, we have\\
\\
\\
$\sum_{k=0}^{n}w_{m,r;\overline{\boldsymbol{\alpha}}}(n,k)(mx+r)^k$
\begin{eqnarray*}
&=&(mx+r)\sum_{k=0}^{n-1}w_{m,r;\overline{\boldsymbol{\alpha}}}(n-1,k)(mx+r)^k\\
&\:\:\:\:\:\:\:\:\:-&(r+m\alpha_{n-1})\sum_{k=0}^{n-1}w_{m,r;\overline{\boldsymbol{\alpha}}}(n-1,k)(mx+r)^k\\
&=&\sum_{k=1}^{n}w_{m,r;\overline{\boldsymbol{\alpha}}}(n-1,k)(mx+r)^k\\
&\:\:\:\:\:\:\:\:-&(r+m\alpha_{n-1})\sum_{k=0}^{n-1}w_{m,r;\overline{\boldsymbol{\alpha}}}(n-1,k)(mx+r)^k.
\end{eqnarray*}
Equating the coefficients of $(mx+r)^k$ on both sides, we get Eq. (\ref{2}).\\
Using Eq. (\ref{2}) it is easy to prove that $w_{m,r;\overline{\boldsymbol{\alpha}}}(n,0)=(-1)^n\prod_{i=1}^{n}(r+m\alpha_{i-1}).$\\
\\
\textbf{Special cases}:
\begin{enumerate}
\item Setting $\alpha_{i}=0,$ for $i=0,1,\cdots,n-1$, hence Eq. (\ref{1}) is reduced to
\begin{equation}
m^nx^n=\sum_{k=0}^{n}w_{m,r;\boldsymbol{0}}(n,k)(mx+r)^k.
\end{equation}
Thus
\begin{eqnarray*}
m^nx^n&=&\sum_{k=0}^{n}w_{m,r;\boldsymbol{0}}(n,k)\sum_{i=0}^{k}\Big(^k_i\Big)(mx)^ir^{k-i}\\
&=&\sum_{i=0}^{n}\Big(\sum_{k=i}^{n}\Big(^k_i\Big)r^{k-i}w_{m,r;\boldsymbol{0}}(n,k)\Big)(mx)^i,
\end{eqnarray*}
hence
\begin{equation}
\sum_{k=i}^{n}\Big(^k_i\Big)r^{k-i}w_{m,r;\boldsymbol{0}}(n,k)=\delta_{ni},
\end{equation}
where $\boldsymbol{0} = (0,0,\cdots,0)$ and $\delta_{ni}$ is Kronecker's delta.
\item Setting $\alpha_i=\alpha,$ for $i=0,1,\cdots,n-1,$ hence Eq. (\ref{1}) is reduced to
\begin{equation}
m^n(x-\alpha)^n=\sum_{k=0}^{n}w_{m,r;\boldsymbol{\alpha}}(n,k)(mx+r)^k,
\end{equation}
therefore we have
\begin{eqnarray*}
m^n\sum_{i=0}^{n}\Big(^n_i\Big)x^i(-\alpha)^{n-i}&=&\sum_{k=0}^{n}w_{m,r;\boldsymbol{\alpha}}(n,k)\sum_{i=0}^{k}\Big(^k_i\Big)(mx)^ir^{k-i}\\
&=&\sum_{i=0}^{n}\Big(\sum_{k=i}^{n}\Big(^k_i\Big)r^{k-i}w_{m,r;\boldsymbol{\alpha}}(n,k)\Big)m^ix^i.
\end{eqnarray*}
Equating the coefficient of $x^i$ on both sides, we get
\begin{equation}
m^{n-i}\Big(^n_i\Big)(-\alpha)^{n-i}=\sum_{k=i}^{n}\Big(^k_i\Big)r^{k-i}w_{m,r;\boldsymbol{\alpha}}(n,k),
\end{equation}
where $\boldsymbol{\alpha}=(\alpha, \alpha,\cdots, \alpha)$.
\item Setting $\alpha_i=i,$ for $i=0,1,\cdots,n-1,$ hence Eq. (\ref{1}) and Eq. (\ref{2}), respectively, are reduced to
\begin{equation}
m^n(x)_n=\sum_{k=0}^{n}w_{m,r;\boldsymbol{i}}(n,k)(mx+r)^k,
\end{equation}
and
\begin{equation}
w_{m,r;\boldsymbol{i}}(n,k)=w_{m,r;\boldsymbol{i}}(n-1,k-1)-(r+m(n-1))w_{m,r;\boldsymbol{i}}(n-1,k),
\end{equation}
where $\boldsymbol{i}=(0,1,\cdots,n-1)$ and $w_{m,r;\boldsymbol{i}}(n,k) = w_{m,r}(n,k)$ are the r-Whitney numbers of the first kind.
\item Setting $\alpha_i=i,$ for $i=0,1,\cdots,n-1,$ and $r=-a$ hence $w_{m,-a;\boldsymbol{\alpha}}$ are the noncentral Whitney number of the first kind, see \cite{Man1}.
\item Setting $\alpha_i=i^p,$ for $i=0,1,\cdots,n-1,$ hence Eq. (\ref{1}) is reduced to
\begin{equation}
m^n\prod_{i=0}^{n}(x-i^p)=\sum_{k=0}^{n}w_{m,r;\boldsymbol{i}^p}(n,k)(mx+r)^k.
\end{equation}
Sun \cite{Sun1} defined p-Stirling numbers of the first kind as
\[
\prod_{i=0}^{n}(x-i^p)=\sum_{i=0}^{n}s_1(n,i,p)x^i,
\]
therefore, we have
\begin{eqnarray*}
m^n\sum_{i=0}^{n}s_1(n,i,p)x^i&=&\sum_{k=0}^{n}w_{m,r;\boldsymbol{i}^p}(n,k)(mx+r)^k\\
&=&\sum_{k=0}^{n}w_{m,r;\boldsymbol{i}^p}(n,k)\sum_{i=0}^{k}\Big(^k_i\Big)(mx)^ir^{k-i}\\
&=&\sum_{i=0}^{n}\Big(\sum_{k=i}^{n}\Big(^k_i\Big)r^{k-i}w_{m,r;\boldsymbol{i}^p}(n,k)\Big)(mx)^i.
\end{eqnarray*}
Equating the coefficient of $x^i$ on both sides, we get
\begin{equation}
m^{n-i}s_1(n,i,p)=\sum_{k=i}^{n}\Big(^k_i\Big)r^{k-i}w_{m,r;\boldsymbol{i}^p}(n,k),
\end{equation}
where $\boldsymbol{i}^p=(0^p,1^p,\cdots,(n-1)^p)$.
\item Setting $\alpha_i=[\alpha_i]_q,$ and $x=[x]_q$ for $i=0,1,\cdots,n-1,$ Eq. (\ref{1}) is reduced to
\begin{equation}
m^n([x]_q;[\overline{\boldsymbol{\alpha}}]_q)_n=\sum_{k=0}^{n}w_{m,r;[\overline{\boldsymbol{\alpha}}]_q}(n,k)(m[x]_q+r)^k,
\end{equation}
where $[\overline{\boldsymbol{\alpha}}]_q=([\alpha_0]_q,[\alpha_1]_q,\cdots,[\alpha_{n-1}]_q)$.\\
El-Desouky and Gomaa \cite{Des3} defined the generalized q-Stirling numbers of the first kind by
\begin{equation}
[x;\overline{\boldsymbol{\alpha}}]_{\underline{n},q}=q^{-\sum_{i=0}^{n-1}\alpha_i}\sum_{k=0}^{n}s_{q,\overline{\boldsymbol{\alpha}}}(n,k)[x]_{q}^{k},
\end{equation}
hence, we get
\[
([x]_q;[\overline{\boldsymbol{\alpha}}]_q)_n=\sum_{k=0}^{n}s_{q,\overline{\boldsymbol{\alpha}}}(n,k)[x]_{q}^{k},
\]
thus we have
\begin{eqnarray*}
m^n\sum_{i=0}^{n}s_{q,\overline{\boldsymbol{\alpha}}}(n,i)[x]_{q}^{i}&=&\sum_{k=0}^{n}w_{m,r;[\overline{\boldsymbol{\alpha}}]_q}(n,k)\sum_{i=0}^{k}\Big(^k_i\Big)r^{k-i}(m[x]_q)^i\\
&=&\sum_{i=0}^{n}\Big(\sum_{k=i}^{n}\Big(^k_i\Big)r^{k-i}w_{m,r;[\overline{\boldsymbol{\alpha}}]_q}(n,k)\Big)m^i[x]_{q}^{i}.
\end{eqnarray*}
Equating the coefficient of $[x]_{q}^{i}$ on both sides, we get
\begin{equation}
m^{n-i}s_{q,\overline{\boldsymbol{\alpha}}}(n,i)=\sum_{k=i}^{n}\Big(^k_i\Big)r^{k-i}w_{m,r;[\overline{\boldsymbol{\alpha}}]_q}(n,k).
\end{equation}
\end{enumerate}
\section{The generalized r-Whitney numbers of the second kind}
\label{sec:2}
\begin{definition} The generalized r-Whitney numbers of the second kind $W_{m,r}(n,k;\overline{\boldsymbol{\alpha}})$ with parameter $\overline{\boldsymbol{\alpha}}=(\alpha_0, \alpha_1,\cdots, \alpha_{n-1})$ are defined by
\begin{equation}\label{3}
(mx+r)^n=\sum_{k=0}^{n}W_{m,r;\overline{\boldsymbol{\alpha}}}(n,k)m^k(x;\overline{\boldsymbol{\alpha}})_k,
\end{equation}
where $W_{m,r;\overline{\boldsymbol{\alpha}}}(0,0)=1$ and $W_{m,r;\overline{\boldsymbol{\alpha}}}(n,k)=0$ for $k>n$.
\end{definition}
\begin{theorem}
The generalized r-Whitney numbers of the second kind $W_{m,r;\overline{\boldsymbol{\alpha}}}(n,k)$ satisfy the recurrence relation
\begin{equation}\label{4}
W_{m,r;\overline{\boldsymbol{\alpha}}}(n,k)=W_{m,r;\overline{\boldsymbol{\alpha}}}(n-1,k-1)+(r+m\alpha_{k})W_{m,r;\overline{\boldsymbol{\alpha}}}(n-1,k)
\end{equation}
for $k\geq1$, and $W_{m,r;\overline{\boldsymbol{\alpha}}}(n,0)=(r+m\alpha_0)^n$.
\end{theorem}
{ \it Proof.}
Since $(mx+r)^n=(mx-m\alpha_k+m\alpha_k+r)(mx+r)^{n-1},$ we have\\
\\
$\sum_{k=0}^{n}W_{m,r;\overline{\boldsymbol{\alpha}}}(n,k)m^k(x;\overline{\boldsymbol{\alpha}})_k$
\begin{eqnarray*}
&=&m(x-\alpha_k)\sum_{k=0}^{n-1}W_{m,r;\overline{\boldsymbol{\alpha}}}(n-1,k)m^k(x;\overline{\boldsymbol{\alpha}})_k\\
&\:\:\:\:\:\:\:\:+&(r+m\alpha_k)\sum_{k=0}^{n-1}W_{m,r;\overline{\boldsymbol{\alpha}}}(n-1,k)m^k(x;\overline{\boldsymbol{\alpha}})_k\\
&=&\sum_{k=1}^{n}W_{m,r;\overline{\boldsymbol{\alpha}}}(n-1,k)m^k(x;\overline{\boldsymbol{\alpha}})_k\\
&\:\:\:\:\:\:\:\:+&(r+m\alpha_k)\sum_{k=0}^{n-1}W_{m,r;\overline{\boldsymbol{\alpha}}}(n-1,k)m^k(x;\overline{\boldsymbol{\alpha}})_k.
\end{eqnarray*}
Equating the coefficient of $(x;\overline{\boldsymbol{\alpha}})_k$ on both sides, we get Eq. (\ref{4}).\\
From Eq. (\ref{4}) it is easy to prove that $W_{m,r;\overline{\boldsymbol{\alpha}}}(n,0)=(r+m\alpha_0)^n$.\\
\begin{theorem}
The generalized r-Whitney numbers of the second kind have the exponential generating function
\begin{equation}\label{5}
\varphi_k(t;\overline{\boldsymbol{\alpha}})=\sum_{i=0}^{k}\frac{e^{(r+m\alpha_i)t}}{\prod_{j=0,j\neq{i}}^{k}(m\alpha_i-m\alpha_j)}.
\end{equation}
\end{theorem}
{ \it Proof.}
The exponential generating function of $W_{m,r;\overline{\boldsymbol{\alpha}}}(n,k)$ is defined by
\begin{equation}\label{6}
\varphi_k(t;\overline{\boldsymbol{\alpha}})=\sum_{n=k}^{\infty}W_{m,r;\overline{\boldsymbol{\alpha}}}(n,k)\frac{t^n}{n!}
\end{equation}
where $W_{m,r;\overline{\boldsymbol{\alpha}}}(n,k)=0$ for $n<k$.
If $k=0$ we have
\[
\varphi_0(t;\overline{\boldsymbol{\alpha}})=\sum_{n=0}^{\infty}W_{m,r;\overline{\boldsymbol{\alpha}}}(n,0)\frac{t^n}{n!}=\sum_{n=0}^{\infty}(r+m\alpha_0)^n\frac{t^n}{n!}=e^{(r+m\alpha_0)t}.
\]
Differentiating both sides of Eq. (\ref{6}) with respect to t. We get
\begin{equation}
\grave{\varphi}_k(t;\overline{\boldsymbol{\alpha}})=\sum_{n=k}^{\infty}W_{m,r;\overline{\boldsymbol{\alpha}}}(n,k)\frac{t^{n-1}}{(n-1)!}
\end{equation}
and from (\ref{4}) we have
\begin{eqnarray*}
\grave{\varphi}_k(t;\overline{\boldsymbol{\alpha}})&=&\sum_{n=k}^{\infty}W_{m,r;\overline{\boldsymbol{\alpha}}}(n-1,k-1)\frac{t^{n-1}}{(n-1)!}\\
&+&(r+m\alpha_k)\sum_{n=k}^{\infty}W_{m,r;\overline{\boldsymbol{\alpha}}}(n-1,k)\frac{t^{n-1}}{(n-1)!}\\
&=&\varphi_{k-1}(t;\overline{\boldsymbol{\alpha}})+(r+m\alpha_k)\varphi_k(t;\overline{\boldsymbol{\alpha}}).
\end{eqnarray*}
Solving this difference-differential equation we obtain Eq. (\ref{5}).\\
\begin{theorem}
The generalized r-Whitney numbers of the second kind have the explicit formula
\begin{equation}\label{19}
W_{m,r;\overline{\boldsymbol{\alpha}}}(n,k)=\sum_{i=0}^k\frac{1}{m^k(\alpha_i)_k}(r+m\alpha_i)^n.
\end{equation}
\end{theorem}
{ \it Proof.}
The exponential generating of $W_{m,r;\overline{\boldsymbol{\alpha}}}(n,k)$ is
\[
\varphi_k(t;\overline{\boldsymbol{\alpha}})=\sum_{i=0}^{k}\frac{e^{(r+m\alpha_i)t}}{\prod_{j=0,j\neq{i}}^{k}(m\alpha_i-m\alpha_j)},
\]
then, we get
\begin{eqnarray*}
\sum_{i=0}^{\infty}W_{m,r;\overline{\boldsymbol{\alpha}}}(n,k)\frac{t^n}{n!}&=&\sum_{i=0}^{k}\frac{1}{m^k(\alpha_i)_k}\sum_{n=0}^{\infty}(r+m\alpha_i)^n\frac{t^n}{n!}\\
&=&\sum_{i=0}^{\infty}\Big(\sum_{i=0}^k\frac{1}{m^k(\alpha_i)_k}(r+m\alpha_i)^n\Big)\frac{t^n}{n!}.
\end{eqnarray*}
Equating the coefficient of $t^n$ on both sides, we get Eq. (\ref{19}).\\
\\
\textbf{Special cases}:
\begin{enumerate}
\item Setting $\alpha_{i}=0,$ for $i=0,1,\cdots,n-1$, hence Eq. (\ref{3}) is reduced to
\begin{equation}
(mx+r)^n=\sum_{k=0}^{n}W_{m,r;\boldsymbol{0}}(n,k)m^kx^k,
\end{equation}
therefore
\begin{equation}
\sum_{k=0}^{n}\Big(^n_k\Big)(mx)^kr^{n-k}=\sum_{k=0}^{n}W_{m,r;\boldsymbol{0}}(n,k)m^kx^k.
\end{equation}
Equating the coefficient of $x^k$ on both sides, we get
\begin{equation}
\Big(^n_k\Big)r^{n-k}=W_{m,r;\boldsymbol{0}}(n,k),
\end{equation}
where $W_{m,r;\boldsymbol{0}}(n,k)$ denotes the generalized Pascal numbers, for more details see \cite{Cal1}, \cite{Sta1}.
\item Setting $\alpha_{i}=\alpha,$ for $i=0,1,\cdots,n-1$, hence Eq. (\ref{3}) is reduced to
\begin{equation}
(mx+r)^n=\sum_{k=0}^{n}W_{m,r;\boldsymbol{\alpha}}(n,k)m^k(x-\alpha)^k,
\end{equation}
so
\begin{eqnarray*}
\sum_{i=0}^{n}\Big(^n_i\Big)(mx)^ir^{n-i}&=&\sum_{k=0}^{n}W_{m,r;\boldsymbol{\alpha}}(n,k)m^k\sum_{i=0}^{k}\Big(^k_i\Big)x^i(-\alpha)^{k-i}\\
&=&\sum_{i=0}^{n}\Big(\sum_{k=i}^{n}\Big(^k_i\Big)W_{m,r;\boldsymbol{\alpha}}(n,k)(\alpha)^{k-i}m^k\Big)x^i.
\end{eqnarray*}
Equating the coefficient of $x^i$ on both sides, we get
\begin{equation}
\Big(^n_i\Big)r^{n-i}m^i=\sum_{k=i}^{n}\Big(^k_i\Big)W_{m,r;\boldsymbol{\alpha}}(n,k)(\alpha)^{k-i}m^k.
\end{equation}
\item Setting $\alpha_{i}=i,$ for $i=0,1,\cdots,n-1$, hence Eq. (\ref{3}) and Eq. (\ref{4}), respectively, are reduced to
\begin{equation}
(mx+r)^n=\sum_{k=0}^{n}W_{m,r;\boldsymbol{i}}(n,k)m^k(x)_k,
\end{equation}
and
\begin{equation}
W_{m,r;\boldsymbol{i}}(n,k)=W_{m,r;\boldsymbol{i}}(n-1,k-1)+(r+mk)W_{m,r;\boldsymbol{i}}(n-1,k),
\end{equation}
where $\boldsymbol{i}=(0,1,\cdots,n-1)$ and $W_{m,r;\boldsymbol{i}}(n,k) = W_{m,r}(n,k)$ are the r-Whitney numbers of the second kind.\\
\\
\textbf{Remark}:\\
Setting $\alpha_i=i,  i=0, 1, \cdots, n-1$ in Eq. (\ref{5}) and using the identity $\prod_{j=0,j\neq{i}}^{k}(i-j)=(-1)^{k-i}(k-i)!j!,$ given by Gould \cite{Gou1}, we obtain the exponential generating function of r-Whitney numbers of the second kind, see \cite{Mez1}, \cite{Che1}.
\item Setting $\alpha_i=i,$ for $i=0,1,\cdots,n-1,$ and $r=-a$ hence Eq. (\ref{3}) is reduced to the noncentral Whitney number of the second kind see \cite{Man1}.
\item Setting $\alpha_{i}=\Big(^i_p\Big),$ for $i=p-1,p,\cdots,n+p-2$, hence Eq. (\ref{3}) is reduced to
\begin{equation}
(mx+r)^n=\sum_{k=0}^{n}W_{m,r;\boldsymbol{p}}(n,k)m^k\prod_{i=p-1}^{p+k-1}(x-\Big(^i_p\Big)).
\end{equation}
Sun \cite{Sun1} defined the p-Stirling numbers of the second kind as
\[
\sum_{k=1}^ns_2(n,k,p)x^i=\prod_{i=p-1}^{p+n-2}(x-\Big(^n_i\Big)),
\]
hence we have
\[
(mx+r)^n=\sum_{k=0}^{n}W_{m,r;\boldsymbol{p}}(n,k)m^k\sum_{i=1}^{k}s_2(k,i,p)x^i
\]
\[
\sum_{i=1}^n\Big(^n_i\Big)r^{n-i}m^ix^i=\sum_{i=1}^{n}\Big(\sum_{k=i-1}^{n}W_{m,r;\boldsymbol{p}}(n,k)m^ks_2(k,i,p)\Big)x^i.
\]
Equating the coefficient of $x^i$ on both sides, we get
\[
\Big(^n_i\Big)r^{n-i}m^i=\sum_{k=i-1}^{n}W_{m,r;\boldsymbol{p}}(n,k)m^ks_2(k,i,p),
\]
where $\boldsymbol{p}=(p-1,p,\cdots,n+p-2)$.
\item Setting $\alpha_{i}=[\alpha_i]_q,$ and $x=[x]_q$ for $i=0,1,\cdots,n-1$, hence Eq. (\ref{3}) is reduced to
\begin{equation}
(m[x]_q+r)^n=\sum_{k=0}^{n}W_{m,r;[\overline{\boldsymbol{\alpha}}]_q}(n,k)m^k([x]_q-[\alpha]_q)_{\underline{k}},
\end{equation}
El-Desouky and Gomaa \cite{Des3} defined the generalized q-Stirling numbers of the second kind as
\begin{eqnarray*}
[x]_q^i&=&\sum_{k=0}^{i}q^{\sum_{l=0}^{k-1}\alpha_l}S_{q,\overline{\boldsymbol{\alpha}}}(i,k)[x;\overline{\boldsymbol{\alpha}}]_{\underline{k},q}\\
&=&\sum_{k=0}^{i}S_{q,\overline{\boldsymbol{\alpha}}}(i,k)([x]_q;[\overline{\boldsymbol{\alpha}}]_q)_{\underline{k}},
\end{eqnarray*}
therefore we have
\begin{eqnarray*}
\sum_{k=0}^{n}W_{m,r;[\overline{\boldsymbol{\alpha}}]_q}(n,k)m^k([x]_q-[\alpha]_q)_{\underline{k}}\\
&=&\sum_{i=0}^{n}\Big(^n_i\Big)r^{n-i}m^i[x]_q^i\\
&=&\sum_{i=0}^{n}\Big(^n_i\Big)r^{n-i}m^i\sum_{k=0}^{i}S_{q,\overline{\boldsymbol{\alpha}}}(i,k)([x]_q;[\overline{\boldsymbol{\alpha}}]_q)_{\underline{k}}\\
&=&\sum_{k=0}^{n}\Big(\sum_{i=k}^{n}\Big(^n_i\Big)r^{n-i}m^iS_{q,\overline{\boldsymbol{\alpha}}}(i,k)\Big)([x]_q;[\overline{\boldsymbol{\alpha}}]_q)_{\underline{k}}.
\end{eqnarray*}
Equating the coefficient of $([x]_q;[\overline{\boldsymbol{\alpha}}]_q)_{\underline{k}}$ on both sides we get
\[
W_{m,r;[\overline{\boldsymbol{\alpha}}]_q}(n,k)m^k=\sum_{i=k}^{n}\Big(^n_i\Big)r^{n-i}m^iS_{q,\overline{\boldsymbol{\alpha}}}(i,k).
\]
\end{enumerate}
\section{The generalized Whitney numbers}
\label{sec:3}
When $r=1$, the generalized r-Whitney numbers of the first $w_{m,1;\overline{\boldsymbol{\alpha}}}(n,k)$ and second kind $W{m,1;\overline{\boldsymbol{\alpha}}}(n,k)$ are reduced, respectively, to numbers which we call the generalized Whitney numbers of the first and second kind which briefly are denoted by $\tilde{w}_{m;\overline{\boldsymbol{\alpha}}}(n,k)$ and $\tilde{W}_{m;\overline{\boldsymbol{\alpha}}}(n,k)$, respectively.
\subsection{The generalized Whitney numbers of the first kind}
\begin{definition} The generalized Whitney numbers of the first kind $\tilde{w}_{m;\overline{\boldsymbol{\alpha}}}(n,k)$ with parameter $\overline{\boldsymbol{\alpha}}=(\alpha_0, \alpha_1,\cdots, \alpha_{n-1})$ are defined by
\begin{equation}\label{18}
m^n(x;\overline{\boldsymbol{\alpha}})_n=\sum_{k=0}^{n}\tilde{w}_{m;\overline{\boldsymbol{\alpha}}}(n,k)(mx+1)^k,
\end{equation}
where $\tilde{w}_{m;\overline{\boldsymbol{\alpha}}}(0,0)=1$ and $\tilde{w}_{m;\overline{\boldsymbol{\alpha}}}(n,k)=0$ for $k>n$.
\end{definition}
\begin{corollary}
The generalized Whitney numbers of the first kind $\tilde{w}_{m;\overline{\boldsymbol{\alpha}}}(n,k)$ satisfy the recurrence relation
\begin{equation}
\tilde{w}_{m;\overline{\boldsymbol{\alpha}}}(n,k)=\tilde{w}_{m;\overline{\boldsymbol{\alpha}}}(n-1,k-1)-(1+m\alpha_{n-1})\tilde{w}_{m;\overline{\boldsymbol{\alpha}}}(n-1,k),
\end{equation}
for $k\geq1$, and $\tilde{w}_{m;\overline{\boldsymbol{\alpha}}}(n,0)=(-1)^n\prod_{i=1}^{n}(1+m\alpha_{i-1})$.
\end{corollary}
{ \it Proof.}
The proof follows directly by setting $r=1$ in Eq. (\ref{2}).\\
\\
\textbf{Special cases}:\\
\begin{enumerate}
\item Setting $\alpha_{i}=0,$ for $i=0,1,\cdots,n-1$  in Eq. (\ref{18}), we get
\begin{eqnarray}
m^nx^n&=&\sum_{k=0}^{n}\tilde{w}_{m;\overline{\boldsymbol{0}}}(n,k)\sum_{i=0}^{k}\Big(^k_i\Big)(mx)^i\\
&=&\sum_{i=0}^{n}\Big(\sum_{k=i}^{n}\Big(^k_i\Big)\tilde{w}_{m;\overline{\boldsymbol{0}}}(n,k)\Big)(mx)^i.
\end{eqnarray}
Hence
\[
\sum_{k=i}^{n}\Big(^k_i\Big)\tilde{w}_{m;\overline{\boldsymbol{0}}}(n,k)=\delta_{ni}.
\]
\item Setting $\alpha_i=\alpha,$  in Eq. (\ref{18}), for $i=0,1,\cdots,n-1,$, we get
\begin{equation}
m^{n-i}\Big(^n_i\Big)(-\alpha)^{n-i}=\sum_{k=i}^{n}\Big(^k_i\Big)\tilde{w}_{m;\boldsymbol{\alpha}}(n,k).
\end{equation}
\item Setting $\alpha_i=i,$ for $i=0,1,\cdots,n-1,$  in Eq. (\ref{18}), we get
\begin{equation}
m^n(x)_n=\sum_{k=0}^{n}\tilde{w}_{m;\boldsymbol{i}}(n,k)(mx+1)^k,
\end{equation}
where $\tilde{w}_{m;\boldsymbol{i}}(n,k) = w_m(n,k)$ are the Whitney numbers of the first kind.
\item Setting $\alpha_i=i^p,$ for $i=0,1,\cdots,n-1,$ in Eq. (\ref{18}), we get
\begin{equation}
m^{n-i}s_1(n,i,p)=\sum_{k=i}^{n}\Big(^k_i\Big)\tilde{w}_{m;\overline{\boldsymbol{i}^p}}(n,k).
\end{equation}
\item Setting $\alpha_i=[\alpha_i]_q,$ and $x=[x]_q$ for $i=0,1,\cdots,n-1,$ in Eq. (\ref{18}), we get
\begin{equation}
m^{n-i}s_{q,\overline{\boldsymbol{\alpha}}}(n,i)=\sum_{k=i}^{n}\tilde{w}_{m;[\overline{\boldsymbol{\alpha}}]_q}(n,k)\Big(^k_i\Big).
\end{equation}
\end{enumerate}
\subsection{The generalized Whitney numbers of the second kind}
\begin{definition} The generalized Whitney numbers of the second kind $\tilde{W}_{m;\overline{\boldsymbol{\alpha}}}(n,k)$ with parameter $\overline{\boldsymbol{\alpha}}=(\alpha_0, \alpha_1,\cdots, \alpha_{n-1})$ are defined by
\begin{equation}\label{20}
(mx+1)^n=\sum_{k=0}^{n}\tilde{W}_{m;\overline{\boldsymbol{\alpha}}}(n,k)m^k(x;\overline{\boldsymbol{\alpha}})_k,
\end{equation}
where $\tilde{W}_{m;\overline{\boldsymbol{\alpha}}}(0,0)=1$ and $\tilde{W}_{m;\overline{\boldsymbol{\alpha}}}(n,k)=0$ for $k>n$.
\end{definition}
\begin{corollary}
The generalized Whitney numbers of the second kind $\tilde{W}_{m;\overline{\boldsymbol{\alpha}}}(n,k)$ satisfy the recurrence relation
\begin{equation}
\tilde{W}_{m;\overline{\boldsymbol{\alpha}}}(n,k)=\tilde{W}_{m;\overline{\boldsymbol{\alpha}}}(n-1,k-1)+(1+m\alpha_{k})\tilde{W}_{m;\overline{\boldsymbol{\alpha}}}(n-1,k),
\end{equation}
for $k\geq1$, and $\tilde{W}_{m;\overline{\boldsymbol{\alpha}}}(n,0)=(1+m\alpha_0)^n$.
\end{corollary}
{ \it Proof.} The proof follows directly by setting $r=1$ in Eq. (\ref{4}).\\
\begin{corollary}
The generalized Whitney numbers of the second kind have the exponential generating function
\begin{equation}
\varphi_k(t;\overline{\boldsymbol{\alpha}})=\sum_{i=0}^{k}\frac{e^{(1+m\alpha_i)t}}{\prod_{j=0,j\neq{i}}^{k}(m\alpha_i-m\alpha_j)}.
\end{equation}
\end{corollary}
{ \it Proof.} The proof follows directly by setting $r=1$ in Eq. (\ref{5}).\\
\\
\textbf{Special cases}:\\
\begin{enumerate}
\item Setting $\alpha_{i}=0,$ for $i=0,1,\cdots,n-1$, in Eq. (\ref{20}), then we get
\begin{equation}
\Big(^n_k\Big)=\tilde{W}_{m;\overline{\boldsymbol{0}}}(n,k),
\end{equation}
where $\tilde{W}_{m;\overline{\boldsymbol{0}}}(n,k)$ are Pascal numbers.
\item Setting $\alpha_{i}=\alpha,$ for $i=0,1,\cdots,n-1$, in Eq. (\ref{20}), then we get
\begin{equation}
\Big(^n_i\Big)m^i=\sum_{k=i}^{n}\Big(^k_i\Big)\tilde{W}_{m;\boldsymbol{\alpha}}(n,k)m^k.
\end{equation}
\item Setting $\alpha_{i}=i,$ for $i=0,1,\cdots,n-1$, in Eq. (\ref{20}), then we get
\begin{equation}
(mx+1)^n=\sum_{k=0}^{n}\tilde{W}_{m;\boldsymbol{i}}(n,k)m^k(x)_k,
\end{equation}
where $\tilde{W}_{m;\boldsymbol{i}}(n,k) = W_m(n,k)$ are the Whitney numbers of the second kind.\\
\\
\textbf{Remark}:\\
Setting $\alpha_i=0$ and $r=1$ in Eq. (\ref{5}) we obtain the exponential  generating function of Whitney numbers of the second kind see \cite{Dow1}.
\item Setting $\alpha_{i}=\Big(^i_p\Big),$ for $i=p-1,p,\cdots,n+p-2$, in Eq. (\ref{20}), then we get
\[
\Big(^n_i\Big)m^i=\sum_{k=i-1}^{n}\tilde{W}_{m;\overline{\boldsymbol{p}}}(n,k)m^ks_2(k,i,p).
\]
\item Setting $\alpha_i=[\alpha_i]_q,$ and $x=[x]_q$ for $i=0,1,\cdots,n-1,$ in Eq. (\ref{20}), then we get
\begin{equation}
\tilde{W}_{m;[\overline{\boldsymbol{\alpha}}]_q}(n,k)m^k=\sum_{i=k}^{n}\Big(^n_i\Big)m^iS_{q,\overline{\boldsymbol{\alpha}}}(i,k).
\end{equation}
\end{enumerate}
\section{Relations between Whitney numbers and some types of numbers}
\label{sec:4}
This section is devoted to drive many important relation between The generalized r-Whitney numbers and different types of Stirling numbers of the first and second kind and the generalized harmonic numbers.
\subsection{The relation between the generalized Stirling numbers of the first kind and the generalized r-Whitney numbers of the first kind}
Comtet \cite{Com1}, \cite{Com2} defined the generalized Stirling numbers of the first kind as
\begin{equation}
(x;\overline{\boldsymbol{\alpha}})_n=\sum_{i=0}^{n}s_{\overline{\boldsymbol{\alpha}}}(n,i)x^i,
\end{equation}
substituting in Eq. (\ref{1}), we obtain
\begin{eqnarray*}
m^n\sum_{i=0}^{n}s_{\overline{\boldsymbol{\alpha}}}(n,i)x^i&=&\sum_{k=0}^{n}w_{m,r;\overline{\boldsymbol{\alpha}}}(n,k)(mx+r)^k\\
&=&\sum_{k=0}^{n}w_{m,r;\overline{\boldsymbol{\alpha}}}(n,k)\sum_{i=0}^{k}\Big(^k_i\Big)(mx)^ir^{k-i}\\
&=&\sum_{i=0}^{n}\Big(\sum_{k=i}^{n}\Big(^k_i\Big)w_{m,r;\overline{\boldsymbol{\alpha}}}(n,k)r^{k-i}\Big)(mx)^i.
\end{eqnarray*}
Equating the coefficient of $x^i$ on both sides, we have
\begin{equation}\label{11}
m^{n-i}s_{\overline{\boldsymbol{\alpha}}}(n,i)=\sum_{k=i}^{n}\Big(^k_i\Big)r^{k-i}w_{m,r;\overline{\boldsymbol{\alpha}}}(n,k).
\end{equation}
Setting $r=1,$ we get
\begin{equation}\label{12}
m^{n-i}s_{\overline{\boldsymbol{\alpha}}}(n,i)=\sum_{k=i}^{n}\Big(^k_i\Big)\tilde{w}_{m;\overline{\boldsymbol{\alpha}}}(n,k).
\end{equation}
\subsection{The relation between the generalized Stirling numbers of the second kind and the generalized r-Whitney numbers of the second kind}
Comtet \cite{Com1} defined the generalized Stirling numbers of the second kind as
\begin{equation}\label{7}
x^k=\sum_{i=0}^{k}S_{\overline{\boldsymbol{\alpha}}}(k,i)(x;\overline{\boldsymbol{\alpha}})_i,
\end{equation}
from Equations (\ref{3}) and (\ref{7}), we have
\begin{eqnarray*}
\sum_{i=0}^{n}W_{m,r;\overline{\boldsymbol{\alpha}}}(n,i)m^i(x;\overline{\boldsymbol{\alpha}})_i&=&\sum_{k=0}^{n}\Big(^n_k\Big)m^kx^kr^{n-k}\\
&=&\sum_{k=0}^{n}\Big(^n_k\Big)m^kr^{n-k}\sum_{i=0}^{k}S_{\overline{\boldsymbol{\alpha}}}(k,i)(x;\overline{\boldsymbol{\alpha}})_i\\
&=&\sum_{i=0}^{n}\Big(\sum_{k=i}^{n}\Big(^n_k\Big)m^kr^{n-k}S_{\overline{\boldsymbol{\alpha}}}(k,i)\Big)(x;\overline{\boldsymbol{\alpha}})_i.
\end{eqnarray*}
Equating the coefficient of $(x;\overline{\boldsymbol{\alpha}})_i$ on both sides, we have
\begin{equation}\label{13}
W_{m,r;\overline{\boldsymbol{\alpha}}}(n,i)m^i=\sum_{k=i}^{n}\Big(^n_k\Big)m^kr^{n-k}S_{\overline{\boldsymbol{\alpha}}}(k,i).
\end{equation}
Setting $r=1,$ we get
\begin{equation}\label{14}
\tilde{W}_{m;\overline{\boldsymbol{\alpha}}}(n,i)m^i=\sum_{k=i}^{n}\Big(^n_k\Big)m^kS_{\overline{\boldsymbol{\alpha}}}(k,i).
\end{equation}
\subsection{The relation between the multiparameter Stirling numbers of the first kind and the generalized r-Whitney numbers of the second kind}
El-Desouky \cite{Des2} defined the multiparameter Stirling numbers of the first kind as
\begin{equation}\label{8}
(x)_n=\sum_{k=0}^{n}s(n,k;\overline{\boldsymbol{\alpha}})(x;\overline{\boldsymbol{\alpha}})_k,
\end{equation}
from Eq. (\ref{3}) we have
\[
(mx+r)^n=\sum_{i=0}^{n}W_{m,r;\overline{\boldsymbol{\alpha}}}(n,i)m^i(x;\overline{\boldsymbol{\alpha}})_i,
\]
hence
\begin{equation}\label{10}
\sum_{k=0}^{n}W_{m,r}(n,k) m^k(x)_k=\sum_{i=0}^{n}W_{m,r;\overline{\boldsymbol{\alpha}}}(n,i)m^i(x;\overline{\boldsymbol{\alpha}})_i,
\end{equation}
from Eq. (\ref{8}) we get
\[
\sum_{k=0}^{n}W_{m,r}(n,k)m^k\sum_{i=0}^{k}s(k,i;\overline{\boldsymbol{\alpha}})(x;\overline{\boldsymbol{\alpha}})_i=\sum_{i=0}^{n}W_{m,r;\overline{\boldsymbol{\alpha}}}(n,i)m^i(x;\overline{\boldsymbol{\alpha}})_i,
\]
then
\[
\sum_{i=0}^{n}\Big(\sum_{k=i}^nm^kW_{m,r}(n,k)s(k,i;\overline{\boldsymbol{\alpha}})\Big)(x;\overline{\boldsymbol{\alpha}})_i=\sum_{i=0}^{n}W_{m,r;\overline{\boldsymbol{\alpha}}}(n,i)m^i(x;\overline{\boldsymbol{\alpha}})_i,
\]
equating the coefficient of $(x;\overline{\boldsymbol{\alpha}})_i$ on both sides, we get
\begin{equation}\label{15}
\sum_{k=i}^nm^{k-i}W_{m,r}(n,k)s(k,i;\overline{\boldsymbol{\alpha}})=W_{m,r;\overline{\boldsymbol{\alpha}}}(n,i).
\end{equation}
Setting $r=1,$ we get
\begin{equation}
\sum_{k=i}^nm^{k-i}W_{m}(n,k)s(k,i;\overline{\boldsymbol{\alpha}})=\tilde{W}_{m;\overline{\boldsymbol{\alpha}}}(n,i).
\end{equation}
\subsection{The relation between the multiparameter Stirling numbers of the second kind and the generalized r-Whitney numbers of the first kind}
El-Desouky \cite{Des2} defined the multiparameter Stirling numbers of the second kind as
 \begin{equation}\label{9}
(x;\overline{\boldsymbol{\alpha}})_n=\sum_{k=0}^{n}S(n,k;\overline{\boldsymbol{\alpha}})(x)_k,
\end{equation}
substituting in Eq. (\ref{1})
\begin{eqnarray*}
m^n\sum_{i=0}^{n}S(n,i;\overline{\boldsymbol{\alpha}})(x)_i&=&\sum_{k=0}^{n}w_{m,r;\overline{\boldsymbol{\alpha}}}(n,k)(mx+r)^k\\
&=&\sum_{k=0}^{n}w_{m,r;\overline{\boldsymbol{\alpha}}}(n,k)\sum_{i=0}^{k}W_{m,r}(k,i) m^i(x)_i\\
&=&\sum_{i=0}^n\Big(\sum_{k=i}^nw_{m,r;\overline{\boldsymbol{\alpha}}}(n,k)W_{m,r}(k,i)\Big)m^i(x)_i.
\end{eqnarray*}
Equating the coefficient of $(x)_i$ on both sides, we get
\begin{equation}\label{16}
m^{n-i}S(n,i;\overline{\boldsymbol{\alpha}})=\sum_{k=i}^nw_{m,r;\overline{\boldsymbol{\alpha}}}(n,k)W_{m,r}(k,i).
\end{equation}
Also, setting $r=1,$ we get
\begin{equation}
m^{n-i}S(n,i;\overline{\boldsymbol{\alpha}})=\sum_{k=i}^n\tilde{w}_{m;\overline{\boldsymbol{\alpha}}}(n,k)W_{m}(k,i).
\end{equation}
\subsection{The relation between the multiparameter Stirling numbers and the generalized r-Whitney numbers of the second kind}
 Similarly, from Eq. (\ref{10}) and Eq. (\ref{9}), we have
\begin{eqnarray*}
\sum_{k=0}^{n}W_{m,r}(n,k)m^k(x)_k&=&\sum_{i=0}^{n}W_{m,r;\overline{\boldsymbol{\alpha}}}(n,i)m^i\sum_{k=0}^iS(i,k;\overline{\boldsymbol{\alpha}})(x)_k\\
&=&\sum_{k=0}^n\Big(\sum_{i=k}^nW_{m,r;\overline{\boldsymbol{\alpha}}}(n,i)S(i,k;\overline{\boldsymbol{\alpha}})m^i\Big)(x)_k.
\end{eqnarray*}
Equating the coefficient of $(x)_k$ on both sides, we get
\begin{equation}\label{17}
W_{m,r}(n,k)=\sum_{i=k}^nW_{m,r;\overline{\boldsymbol{\alpha}}}(n,i)S(i,k;\overline{\boldsymbol{\alpha}})m^{i-k}.
\end{equation}
Setting $r=1,$ we get
\begin{equation}
W_{m}(n,k)=\sum_{i=k}^n\tilde{W}_{m;\overline{\boldsymbol{\alpha}}}(n,i)S(i,k;\overline{\boldsymbol{\alpha}})m^{i-k}.
\end{equation}
\subsection{The relation between the generalized r-Whitney numbers and the generalized Harmonic numbers}
Caki\'{c} \cite{Cak1} defined the generalized harmonic numbers as
\[
H_n(k;\overline{\boldsymbol{\alpha}})=\sum_{i=0}^{n-1}\frac{1}{(\alpha_i)^k}.
\]
From Eq (\ref{1}), we have
\begin{eqnarray}\label{22}
m^n(x;\overline{\boldsymbol{\alpha}})_n&=&\sum_{k=0}^{n}w_{m,r;\overline{\boldsymbol{\alpha}}}(n,k)(mx+r)^k \nonumber\\
&=&\sum_{k=0}^{n}w_{m,r;\overline{\boldsymbol{\alpha}}}(n,k)\sum_{j=0}^{k}\binom{k}{j}(mx)^jr^{k-j}\nonumber\\
&=&\sum_{j=0}^{n}\Big(\sum_{k=j}^{n}\binom{k}{j}r^{k-j}w_{m,r;\overline{\boldsymbol{\alpha}}}(n,k)\Big)(mx)^j,\nonumber \\
&=&\sum_{j=0}^{\infty}\Big(\sum_{k=j}^{n}\binom{k}{j}r^{k-j}w_{m,r;\overline{\boldsymbol{\alpha}}}(n,k)\Big)(mx)^j.
\end{eqnarray}
Also,
\begin{eqnarray*}
m^n(x;\overline{\boldsymbol{\alpha}})_n&=&m^n(x-\alpha_0)(x-\alpha_1) \cdots (x-\alpha_{n-1})\\
&=&m^n\prod_{i=0}^{n}(-\alpha_i)(1-\frac{x}{\alpha_i})\\
&=&m^n\prod_{i=0}^{n}(-\alpha_i).\exp\Big(\sum_{i=0}^{n-1}\log(1-\frac{x}{\alpha_i})\Big)\\
&=&m^n\prod_{i=0}^{n}(-\alpha_i).\exp\Big(-\sum_{i=0}^{n-1}\sum_{k=1}^{\infty}\frac{1}{k}\Big(\frac{x}{\alpha_i}\Big)^k\Big)\\
&=&m^n\prod_{i=0}^{n}(-\alpha_i).\exp\Big(-\sum_{k=1}^{\infty}\frac{x^k}{k}\sum_{i=0}^{n-1}\Big(\frac{1}{\alpha_i}\Big)^k\Big)\\
&=&m^n\prod_{i=0}^{n}(-\alpha_i).\exp\Big(-\sum_{k=1}^{\infty}\frac{x^k}{k}H_n(k;\overline{\boldsymbol{\alpha}})\Big)\\
&=&m^n\prod_{i=0}^{n}(-\alpha_i)\sum_{\ell=0}^{\infty}\frac{(-1)^{\ell}}{\ell!}\Big(\sum_{k=1}^{\infty}\frac{x^k}{k}H_n(k;\overline{\boldsymbol{\alpha}})\Big)^{\ell},\\
\end{eqnarray*}
using Cauchy rule product, this lead to
\begin{eqnarray*}
\Big(\sum_{k=1}^{\infty}\frac{x^k}{k}H_=n(k;\overline{\boldsymbol{\alpha}})\Big)^{\ell}&=&\prod_{j=1}^{\ell}\Big(\sum_{k_j=1}^{\infty}\frac{x^{k_j}}{k_j}H_n(k_j;\overline{\boldsymbol{\alpha}})\Big)\\
&=&\sum_{j=\ell}^{\infty}\Big(\sum_{k_1+k_2+ \cdots +k_{\ell}=j}\frac{1}{k_{1}k_{2} \cdots k_{\ell}}\prod_{\iota=1}^{\ell}H_n(k_\iota;\overline{\boldsymbol{\alpha}})\Big)x^j,\\
\end{eqnarray*}
therefore
\begin{eqnarray}\label{23}
m^n(x;\overline{\boldsymbol{\alpha}})_n&=&m^n\prod_{i=0}^{n}(-\alpha_i)\sum_{\ell=0}^{\infty}\frac{(-1)^{\ell}}{\ell!}\sum_{j=\ell}^{\infty}\Big(\sum_{k_1+k_2+ \cdots +k_{\ell}=j}\frac{1}{k_{1}k_{2} \cdots k_{\ell}}\prod_{\iota=1}^{\ell}H_n(k_\iota;\overline{\boldsymbol{\alpha}})\Big)x^j \nonumber\\
&=&m^n\prod_{i=0}^{n}(-\alpha_i)\sum_{j=0}^{\infty}\Big(\sum_{\ell=0}^{\infty}\frac{(-1)^{\ell}}{\ell!}\Big(\sum_{k_1+k_2+\cdots +k_{\ell}=j}\frac{1}{k_{1}k_{2} \cdots k_{\ell}}\prod_{\iota=1}^{\ell}H_n(k_\iota;\overline{\boldsymbol{\alpha}})\Big)x^j.
\end{eqnarray}
From Eq. (\ref{22}) and Eq. (\ref{23}) we have the following identity
\begin{equation}\label{21}
\sum_{k=j}^{n}\binom{k}{j}r^{k-j}w_{m,r;\overline{\boldsymbol{\alpha}}}(n,k)=m^{n-j}\prod_{i=0}^{n}(-\alpha_i)\sum_{\ell=0}^{\infty}\frac{(-1)^{\ell}}{\ell!}\Big(\sum_{k_1+k_2+ \cdots +k_{\ell}=j}\frac{1}{k_{1}k_{2} \cdots k_{\ell}}\prod_{\iota=1}^{\ell}H_n(k_\iota;\overline{\boldsymbol{\alpha}})\Big),
\end{equation}
From Eq. (\ref{11}) and Eq. (\ref{21}) we have
\begin{equation}
s_{\overline{\boldsymbol{\alpha}}}(n,j)=\prod_{i=0}^{n}(-\alpha_i)\sum_{\ell=0}^{\infty}\frac{(-1)^{\ell}}{\ell!}\Big(\sum_{k_1+k_2+ \cdots +k_{\ell}=j}\frac{1}{k_{1}k_{2} \cdots k_{\ell}}\prod_{\iota=1}^{\ell}H_n(k_\iota;\overline{\boldsymbol{\alpha}})\Big),
\end{equation}
this equation gives the generalized Stirling numbers in terms of the generalized Harmonic numbers.
\section{Matrix representation}
\label{sec:6}
In this section we drive a matrix representation for some given relations.
\begin{enumerate}
\item Eq. (\ref{15}) can be represented in matrix form as
\begin{equation}\label{24}
\widehat{\boldsymbol{W}}_{m,r}\:\boldsymbol{s}(\overline{\alpha})=\widehat{\boldsymbol{W}}_{m,r;\overline{\boldsymbol{\alpha}}},
\end{equation}
where $\widehat{W}_{m,r}(n,k)=m^k W_{m,r}(n,k)$ and $\widehat{W}_{m,r;\overline{\boldsymbol{\alpha}}}(n,i)=m^i W_{m,r;\overline{\boldsymbol{\alpha}}}(n,i)$ and $\boldsymbol{W}_{m,r},\boldsymbol{s}(\overline{\boldsymbol{\alpha}})$ and $\boldsymbol{W}_{m,r;\overline{\boldsymbol{\alpha}}}$ are $n\times n$ lower triangle matrices whose entries are, respectively, the r-Whitney numbers of the second kind, the multiparameter Stirling numbers of the first kind and the generalized r-Whitney numbers of the second kind.\\
For example if $0\leq n,k,i \leq 3,$ and using matrix representation given in \cite{Cak2}, hence Eq. (\ref{24}) can be written as
\[\scriptsize
\left(
\begin{array}{llll}
1&0&0&0\\
r&m&0&0\\
r^2&m(2r+m)&m^2&0\\
r^3&m(3r^2+3mr+m^2)&m^2(3r+3m)&m^3
\end{array}
\right)
\left(
\begin{array}{rrrr}
1&0&0&0\\
\alpha_0&1&0&0\\
\alpha_0(\alpha_0-1)&\alpha_0+\alpha_1-1&1&0\\
s(3,0;\overline{\boldsymbol{\alpha}})&s(3,1;\overline{\boldsymbol{\alpha}})&s(3,2;\overline{\boldsymbol{\alpha}})&1
\end{array}
\right)
\]
\[\scriptsize
\left(
\begin{array}{cccc}
1&0&0&0\\
(r+m\alpha_0)&m&0&0\\
(r+m\alpha_0)^2&(2r+m\alpha_0+m\alpha_1)m&m^2&0\\
(r+m\alpha_0)^3&\widehat{W}_{m,r;\overline{\boldsymbol{\alpha}}}(3,1)&\widehat{W}_{m,r;\overline{\boldsymbol{\alpha}}}(3,2)&m^3
\end{array}
\right)
\]
where $s(3,0;\overline{\boldsymbol{\alpha}})=\alpha_0(\alpha_0-1)(\alpha_0-2),$ $s(3,1;\overline{\boldsymbol{\alpha}})=\alpha_0(\alpha_0-1)+(\alpha_1-2)(\alpha_0+\alpha_1-1),$ $s(3,2;\overline{\boldsymbol{\alpha}})=\alpha_0+\alpha_1+\alpha_2-3,$ $\widehat{W}_{m,r;\overline{\boldsymbol{\alpha}}}(3,1)=((r+m\alpha_0)^2+(r+m\alpha_1)(2r+m\alpha_0+m\alpha_1))m,$ $\widehat{W}_{m,r;\overline{\boldsymbol{\alpha}}}(3,2)=(3r+m\alpha_0+m\alpha_1+m\alpha_2)m^2.$
\item Eq. (\ref{16}) can be represented in a matrix form as
\begin{equation}\label{25}
\boldsymbol{w}_{m,r;\overline{\boldsymbol{\alpha}}}\:\widehat{\boldsymbol{W}}_{m,r}=\widehat{\boldsymbol{S}}(\overline{\boldsymbol{\alpha}}),
\end{equation}
where $\widehat{S}(n,i;\overline{\boldsymbol{\alpha}})=m^n S(n,i;\overline{\boldsymbol{\alpha}}),$ and $\boldsymbol{w}_{m,r;\overline{\boldsymbol{\alpha}}}$ and $\boldsymbol{S}(\overline{\boldsymbol{\alpha}})$ are $n\times n$ lower triangle matrices whose entries are, respectively, the generalized r-Whitney numbers of the first kind and the multiparameter Stirling numbers of the second kind.\\
For example if $0\leq n,k,i \leq 3,$ hence Eq. (\ref{25}) can be written as
\[\scriptsize
\left(
\begin{array}{cccc}
1&0&0&0\\
-r-m\alpha_0&1&0&0\\
(r+m\alpha_0)(r+m\alpha_1)&-2r-m\alpha_0-m\alpha_1&1&0\\
w_{m,r;\overline{\boldsymbol{\alpha}}}(3,0)&w_{m,r;\overline{\boldsymbol{\alpha}}}(3,1)&w_{m,r;\overline{\boldsymbol{\alpha}}}(3,2)&1
\end{array}
\right)
\]
\[\scriptsize
.\left(
\begin{array}{cccc}
1&0&0&0\\
r&m&0&0\\
r^2&m(2r+m)&m^2&0\\
r^3&m(3r^2+3mr+m^2)&m^2(3r+3m)&m^3
\end{array}
\right)
\]
\[\scriptsize
=\left(
\begin{array}{cccc}
1&0&0&0\\
-m\alpha_0&m&0&0\\
m^2\alpha_0\alpha_1&m^2(-\alpha_0-\alpha_1+1)&m^2&0\\
m^3\alpha_0\alpha_1\alpha_2&\widehat{S}(3,1;\overline{\boldsymbol{\alpha}})&\widehat{S}(3,2;\overline{\boldsymbol{\alpha}})&m^3
\end{array}
\right)
\]
where $w_{m,r;\overline{\boldsymbol{\alpha}}}(3,0)=-(r+m\alpha_0)(r+m\alpha_1)(r+m\alpha_2)$, $w_{m,r;\overline{\boldsymbol{\alpha}}}(3,1)=(r+m\alpha_0)(r+m\alpha_1)+(2r+m\alpha_0+m\alpha_1)(r+m\alpha_2)$, $w_{m,r;\overline{\boldsymbol{\alpha}}}(3,2)=-3r-m\alpha_0-m\alpha_1-m\alpha_2,$ $\widehat{S}(3,1;\overline{\boldsymbol{\alpha}})=m^3(\alpha_0\alpha_1+\alpha_0\alpha_2+\alpha_1\alpha_2-\alpha_0-\alpha_1-\alpha_2+1),$ $\widehat{S}(3,2;\overline{\boldsymbol{\alpha}})=m^3(-\alpha_0-\alpha_1-\alpha_2+3).$
\item Eq. (\ref{17}) can be represented in a matrix form as
\begin{equation}\label{26}
\widehat{\boldsymbol{W}}_{m,r;\overline{\boldsymbol{\alpha}}}\:\boldsymbol{S}(\overline{\boldsymbol{\alpha}})=\widehat{\boldsymbol{W}}_{m,r},
\end{equation}
For example if $0\leq n,k,i \leq 3,$ hence Eq. (\ref{25}) can be written as
\[\scriptsize
\left(
\begin{array}{cccc}
1&0&0&0\\
(r+m\alpha_0)&m&0&0\\
(r+m\alpha_0)^2&(2r+m\alpha_0+m\alpha_1)m&m^2&0\\
(r+m\alpha_0)^3&\widehat{W}_{m,r;\overline{\boldsymbol{\alpha}}}(3,1)&\widehat{W}_{m,r;\overline{\boldsymbol{\alpha}}}(3,2)&m^3
\end{array}
\right)
\]
\[\scriptsize
.\left(
\begin{array}{llll}
1&0&0&0\\
-\alpha_0&1&0&0\\
\alpha_0\alpha_1&-\alpha_0-\alpha_1+1&1&0\\
\alpha_0\alpha_1\alpha_2&S(3,1;\overline{\boldsymbol{\alpha}})&S(3,2;\overline{\boldsymbol{\alpha}})&1
\end{array}
\right)
=\left(
\begin{array}{rrrr}
1&0&0&0\\
r&m&0&0\\
r^2&m(2r+m)&m^2&0\\
r^3&m(3r^2+3mr+m^2)&m^2(3r+3m)&m^3
\end{array}
\right)
\]
where $\widehat{W}_{m,r;\overline{\boldsymbol{\alpha}}}(3,1)=((r+m\alpha_0)^2+(r+m\alpha_1)(2r+m\alpha_0+m\alpha_1))m,$ $\widehat{W}_{m,r;\overline{\boldsymbol{\alpha}}}(3,2)=(3r+m\alpha_0+m\alpha_1+m\alpha_2)m^2,$
$S(3,1;\overline{\boldsymbol{\alpha}})=\alpha_0\alpha_1+\alpha_0\alpha_2+\alpha_1\alpha_2-\alpha_0-\alpha_1-\alpha_2+1,$ $S(3,2;\overline{\boldsymbol{\alpha}})=-\alpha_0-\alpha_1-\alpha_2+3.$
\end{enumerate}


\begin{thebibliography}{99}
\bibitem{Ben1}
Benoumhani, M.:
On Whitney numbers of Dowlling lattices.
 Discrete Math. 159, 13--33 (1996)
\bibitem{Cak1}
Caki\'{c},Nenad P. :
The complete Bell polynomials and numbers of Mitrinović.
Univ. Beograd. Publ. Elektrotehn. Fak. Ser. Mat. No 6, 74--78 (1995)
\bibitem{Cak2}
Caki\'{c}, Nenad P., El-Desouky, Beih S., Milovanovi\'{c}, Gradimir V.:
Explicit formulas and combinatorial identities for generalized Stirling numbers.
Mediterr. j. Math., 10 No. 1, 375--385 (2013)
\bibitem{Cal1}
Call, G.S., Vellman, D.J.:
Pascal matrices.
Amer. Math. Monthly 100, 372-–376 (1993)
\bibitem{Cha1}
Charalambides, Ch. A.:
On the q-differences of the generalized q-factorials.
J. Stat. Plan. Infer., 54, 31-34 (1996).
\bibitem{Che1}
Cheon, Gi-Sang, Jung, Ji-Hwan
r-Whitney numbers of Dowling lattices.
Discrete Math. 312, 2337--2348 (2012)
\bibitem{Com1}
Comtet, M.:
Numbers de Stirling generaux et fonctions symetriques.
C. R. Acad. Sc. Par. (Series A) 275, 747--750 (1972)
 \bibitem{Com2}
Comtet, L.:
Advanced Combinatorics: The Art of Finite and Infinite Expansions.
D. Reidel Publishing Company. Dordrecht, Holand (1974)
\bibitem{Dow1}
Dowling, T. A.:
A class of geometric lattices passed on finite groups.
J. Combin. Theory, Ser. B 14, 61--86 (1973)
\bibitem{Des1}
El-Desouky, B. S., Caki\'{c}, Nenad P.
Generalized higher order Stirling numbers.
Math. Comput. Model. 54, 2848–-2857 (2011)
\bibitem{Des2}
El-Desouky, B. S.
The multiparameter non-central Stirling numbers.
Fibonacci Quart 32 (3), 218--225 (1994)
\bibitem{Des3}
El-Desouky, B. S., Gomaa, R. S.
q-Comtet and Generalized q-harmonic Numbers.
Journal of Mathematical Sciences: Advances and Applications, 211, 52--71 (2011)
\bibitem{Gou1}
Gould, H. W.
Combinatorial numbers and associated identities.
table 1: Stirling numbers, Unpublished Manuscript, http://www.math.wvu.edu/gould/vol.7.
\bibitem{Man1}
Mangontarum, Mahid M., Cauntongan, Omar I., Macodi-Ringia, Amila P.
The Noncentral Version of the Whitney Numbers.
International Jornal of Mathematics and Mathematical Science
\bibitem{Mez1}
Mez\H o, I.
A new formula for the Bernoulli polynomials.
Results Math. 58, 329--335 (2010).
\bibitem{Sta1}
Stanimirovi\'{c}, Stefan
A Generalization of The Pascal Matrix and its Properties.
Ser. Math. Inform. 26, 17--27 (2011)
\bibitem{Sun1}
Yidong Sun, "Two Classes of p-Stirling Numbers." Discrete Math., 306(2006) 2801-2805.
\end{thebibliography}
\end{document}